\documentclass[a4paper]{article}
\usepackage{fullpage}
\usepackage{amsmath}
\usepackage{amssymb}
\usepackage{amsthm}
\usepackage{graphicx}
\usepackage{bbm}
\usepackage{enumerate}
\usepackage{microtype}
\usepackage{tocbibind}
\usepackage{subcaption}
\usepackage{xcolor}
\usepackage[none]{hyphenat}
\usepackage[colorlinks=true,urlcolor=blue,linkcolor=blue,plainpages=false,pdfpagelabels]{hyperref}
\allowdisplaybreaks

\newtheorem{theorem}{Theorem}[section]
\newtheorem{proposition}[theorem]{Proposition}
\newtheorem{lemma}[theorem]{Lemma}
\newtheorem{corollary}[theorem]{Corollary}

\newtheorem{notation}[theorem]{Notation}
\newtheorem{definition}[theorem]{Definition}

\newcommand{\N}{\mathbb{N}}

\newcommand{\R}{\mathbb{R}}
\newcommand{\Hh}{\mathbb{H}}
\newcommand{\eps}{\varepsilon}
\newcommand{\vp}{\varphi}

\newcommand{\wt}[1]{\widetilde{#1}}
\newcommand{\vv}[1]{\overrightarrow{#1}}

\newcommand{\prv}[2]{\langle\vv{#1},\vv{#2}\rangle}

\DeclareMathOperator{\Fix}{Fix}
\DeclareMathOperator{\arcosh}{arcosh}

\title{Abstract strongly convergent variants of the proximal point algorithm}
\author{Andrei Sipo\c s${}^{a,b}$\\[2mm]
\footnotesize ${}^a$Research Center for Logic, Optimization and Security (LOS), Department of Computer Science,\\
\footnotesize Faculty of Mathematics and Computer Science, University of Bucharest,\\
\footnotesize Academiei 14, 010014 Bucharest, Romania\\[1mm]
\footnotesize ${}^b$Simion Stoilow Institute of Mathematics of the Romanian Academy,\\
\footnotesize Calea Grivi\c tei 21, 010702 Bucharest, Romania\\[2mm]
\footnotesize E-mail: andrei.sipos@fmi.unibuc.ro\\
}
\date{}

\begin{document}

\maketitle

\begin{abstract}
We prove an abstract form of the strong convergence of the Halpern-type and Tikhonov-type proximal point algorithms in CAT(0) spaces. In addition, we derive uniform and computable rates of metastability (in the sense of Tao) for these iterations using proof mining techniques.


\end{abstract}

\section{Introduction}

The proximal point algorithm is a fundamental tool of convex optimization, usually attributed to Martinet \cite{Mar70}, Rockafellar \cite{Roc76} (who named it) and Br\'ezis and Lions \cite{BreLio78}.  In its many variants, it usually operates by iterating on a starting point -- in, say, a Hilbert space -- a sequence of mappings dubbed ``resolvents'', whose fixed points coincide with the solutions of the optimization problem that one is aiming at. Thus, if for any $\gamma>0$ one denotes the resolvent of order $\gamma$ corresponding to the given problem by $J_\gamma$, then one selects a sequence $(\gamma_n)$ of `step-sizes' and then forms the iterative sequence which bears the name `proximal point algorithm' by putting, for any $n$, $x_{n+1}$ to be equal to $J_{\gamma_n}x_n$.

Unfortunately, this class of algorithms is usually only weakly convergent: that strong convergence does not always hold has been shown by G\"uler \cite{Gul91}. A natural question, then, is how to modify the algorithm into a strongly convergent one. A source of inspiration was found in the iterations commonly used in metric fixed point theory, for example the iteration introduced by Halpern in \cite{Hal67}, which bears his name and which is used to find fixed points of e.g. a self-mapping $T$ of the space and which, for a given `anchor' point $u$ in the space and a sequence of `weights' $(\alpha_n)$, constructs, for each $n$, the step $x_{n+1}$ as
$$\alpha_nu + (1-\alpha_n)Tx_n.$$
In order to guarantee strong convergence of the algorithm, one usually imposes some conditions on the sequence $(\alpha_n)$, for example
$$\lim_{n \to \infty} \alpha_n=0, \quad \sum_{n=0}^\infty \alpha_n = \infty.$$
It was more or less known since Halpern (see \cite[Theorem 6]{Suz09} for an updated proof) that the above two conditions are necessary for strong convergence when $T$ is a nonexpansive mapping, but they may not be sufficient: Halpern himself in his original paper \cite{Hal67} proved strong convergence using some highly restrictive additional conditions which excluded the natural choice $\alpha_n:=1/(n+1)$. Only in the 1990s, Wittmann \cite{Wit92} managed to show strong convergence under a weaker additional condition, which included that choice.

Given this, it is then natural to consider what Kamimura and Takahashi \cite{KamTak00} and Xu \cite{Xu02} independently introduced as the {\it Halpern-type proximal point algorithm}, where the map $T$ in the Halpern iteration above is replaced by the resolvent $J_{\gamma_n}$ from the proximal point algorithm, and which strongly converges in Hilbert spaces if one imposes the Halpern conditions above and the condition $\lim_{n \to \infty} \gamma_n = \infty$. (A related modification known as the {\it Tikhonov regularization} was studied in \cite{LehMou96} and especially in \cite{Xu06}.) Aoyama and Toyoda \cite{AoyToy17} have recently shown that the Halpern proximal point algorithm converges in Banach spaces which are uniformly convex and uniformly smooth if one imposes in addition to the two Halpern conditions above just the condition that the sequence $(\gamma_n)$ is bounded below away from $0$, their proof making use of the property of the resolvents being {\bf strongly nonexpansive}. (In particular, the usual Halpern iteration had already been shown by Saejung \cite{Sae10-sne} to strongly converge for strongly nonexpansive $T$ with just the two Halpern conditions.)

In the last two decades, there has been a continued interest in extending results in fixed point theory and convex optimization from linear spaces like Hilbert or Banach spaces to {\bf nonlinear} ones, chief among them being CAT(0) spaces (to be defined in the next section), frequently regarded as the rightful nonlinear generalization of Hilbert spaces. The first adaptation of the proximal point algorithm to this context was obtained in \cite{Bac13} by Ba\v{c}\'ak, who also authored the book \cite{Bac14}, which serves as a general reference for convex optimization in CAT(0) spaces.

It has been observed by Eckstein that the arguments used to prove the convergence of the usual proximal point algorithm ``hinge primarily on the firmly nonexpansive properties of the resolvents'' \cite[p. 7]{Eck89}. Inspired by this remark, the author, together with Leu\c stean and Nicolae, has introduced in \cite{LeuNicSip18} the concept of {\it jointly firmly nonexpansive family of mappings}, in the context of CAT(0) spaces, which allows for a highly abstract proof of the proximal point algorithm's convergence, encompassing virtually all known variants in the literature. (It was not accidental that we have already presented resolvents above in a quite abstract way.) We have recently revisited \cite{Sip21} the concept, providing a conceptual characterization of it and showing how it may be used to prove other kinds of results usually associated with resolvent-type mappings.

{\it The goal of this paper is to present strongly convergent variants of the proximal point algorithm in the framework of jointly firmly nonexpansive families of mappings in CAT(0) spaces.}

We chose to adapt the proof of Aoyama and Toyoda \cite{AoyToy17} due to the fact that it uses the weakest conditions known so far; even though full strong nonexpansiveness is not available in nonlinear spaces like CAT(0) spaces, what one needs for their proof to go through is the uniform strengthening of the weaker notion of strong quasi-nonexpansiveness, a strengthening which -- with a quantitative modulus -- was introduced by Kohlenbach in \cite{Koh16}.

The quantitative nature of this notion is due to the fact that the investigations in \cite{Koh16} tie into the area of {\it proof mining}, an applied subfield of mathematical logic which aims to analyze proofs in concrete mathematics using tools from proof theory in order to extract additional information from them. This program in its current form has been developed in the last decades primarily by Kohlenbach and his collaborators -- see \cite{Koh08} for a comprehensive monograph; a recent survey which serves as a short and accessible introduction is \cite{Koh19}. It would be natural, then, to ask for a {\it a rate of convergence} for the iterations we mentioned above; unfortunately, rates of convergence for iterative sequences which are commonly employed in nonlinear analysis and convex optimization may not be uniform or computable (see \cite{Neu15}). Kohlenbach's work then suggests that one should look instead at the following (classically but not constructively) equivalent form of the Cauchy property (actually identifiable in mathematical logic as its Herbrand normal form):
$$\forall \eps>0 \,\forall g:\N\to\N \,\exists N\in \N \,\forall i,j\in [N,N+g(N)] \ \left(\| x_i-x_j\| \leq\eps\right),$$
which has been arrived at independently by Tao in his own work on ergodic theory \cite{Tao08A} and popularized in \cite{Tao08} -- as a result of the latter, the property got its name of {\it metastability} (at the suggestion of Jennifer Chayes). Kohlenbach's metatheorems then guarantee the existence of a computable and uniform {\it rate of metastability} -- a bound $\Theta(\eps,g)$ on the $N$ in the sentence above; and this research program of proof mining has achieved over the years a number of non-trivial extractions of such rates from celebrated strong convergence proofs, see, e.g., \cite{Koh11,KohLeu12,KohSip21}.

Some linear space variants of the Halpern-type proximal point algorithm have already been analyzed from the point of view of proof mining -- specifically by Kohlenbach \cite{Koh20} (whose analysis we shall follow closely, given that he analyzed the original proof of Aoyama and Toyoda \cite{AoyToy17}), as well as by Pinto \cite{Pin21} (who analyzed the proof of Xu \cite{Xu02} mentioned above) and Leu\c stean and Pinto \cite{LeuPin21}.

The main obstacle in producing our analysis is, as suggested before, strong nonexpansiveness. We have said that for the usual (non-quantitative) proof we can use only the uniform version of the strong quasi-nonexpansive property -- witnessed here by Lemma~\ref{sne1} -- but this turns out not to be enough for the quantitative version. What we do is to mine the proof of that lemma in order to obtain a further `quantitative quasiness' property in the form of Proposition~\ref{sne-q}, which gives us exactly the necessary ingredient for the proof to go through, namely the analogue of \cite[Lemma 3.3]{Koh20} in Kohlenbach's original analysis.

Section~\ref{sec:prelim} presents the general concepts we shall need regarding CAT(0) spaces, self-mappings of them (including the jointly firmly nonexpansive families mentioned above) and techniques to prove convergence. We chose to present the qualitative convergence results distinctly from the quantitative ones, so that they could stand on their own. Thus, the main convergence theorems can be found in Section~\ref{sec:main} -- Theorem~\ref{main} and Corollary~\ref{main-t}, showing the strong convergence of the Halpern-type and of the Tikhonov-type proximal point algorithm, respectively -- while their corresponding quantitative versions, yielding rates of metastability, can be found in Section~\ref{sec:quant}.

\section{Preliminaries}\label{sec:prelim}

One says that a metric space $(X,d)$ is {\it geodesic} if for any two points $x$, $y \in X$ there is a {\it geodesic} that joins them, i.e. a mapping $\gamma : [0,1] \to X$ such that $\gamma(0)=x$, $\gamma(1)=y$ and for any $t$, $t' \in [0,1]$ we have that
$$d(\gamma(t),\gamma(t')) = |t-t'| d(x,y).$$
Among geodesic spaces, a subclass that is usually considered (e.g. in convex optimization) to be the rightful nonlinear generalization of Hilbert spaces is the class of CAT(0) spaces, introduced by A. Aleksandrov \cite{Ale51} and named as such by M. Gromov \cite{Gro87}, defined as those geodesic spaces $(X,d)$ such that for any geodesic $\gamma : [0,1] \to X$ and for any $z \in X$ and $t \in [0,1]$ we have that
$$d^2(z,\gamma(t)) \leq (1-t)d^2(z,\gamma(0)) + td^2(z,\gamma(1)) - t(1-t)d^2(\gamma(0),\gamma(1)).$$
Another well-known fact about CAT(0) spaces is that each such space $(X,d)$ is {\it uniquely geodesic} -- that is, for any $x$, $y \in X$ there is a unique  geodesic $\gamma : [0,1] \to X$ that joins them -- and in this context we shall denote, for any $t \in [0,1]$, the point $\gamma(t)$ by $(1-t)x + ty$. Note also that any CAT(0) space $(X,d)$ is {\it Busemann convex} -- i.e., for any $x$, $y$, $u$, $v \in X$ and $t \in [0,1]$, 
$$d((1-t)x + ty,(1-t)u + tv) \leq (1-t)d(x,u) + td(y,v),$$
and in particular, for any $x$, $u$, $v \in X$ and $t \in [0,1]$, 
$$d(x,(1-t)u + tv) \leq (1-t)d(x,u) + td(x,v).$$

The simplest example of a non-Hilbert (complete) CAT(0) space is the {\it (hyperbolic) Poincar\'e upper half-plane model} (see also \cite{BriHae99}), having as the underlying set
$$\Hh:= \{(x,y) \in \R^2 \mid y > 0\},$$
where, given the function $\arcosh : [1,\infty) \to [0,\infty)$, where for every $x \in [1,\infty)$, $\arcosh x = \ln (x+\sqrt{x^2-1})$, the distance function is defined as follows: for any $(x_1,y_1)$, $(x_2,y_2) \in \Hh$, one sets
$$d((x_1,y_1),(x_2,y_2)):= \arcosh \left( 1 + \frac{(x_2-x_1)^2 + (y_2-y_1)^2}{2y_1y_2} \right).$$
One proves, usually by considering this distance function as arising from a Riemannian metric, that geodesic lines of this space are of two types: for every $a\in \R$ and $r > 0$, one has the semicircle
$$\mathcal{C}_{a,r} = \{(x,y) \in \Hh \mid (x-a)^2+y^2=r^2\},$$
while for every $a \in \R$, one has the ray
$$\mathcal{R}_a = \{(x,y) \in \Hh \mid x=a\}.$$
It may be then easily shown that for every two points there is exactly one geodesic segment that joins them. The exact formulas for the convex combination of two points (as well as those for projecting onto geodesics) are somewhat involved and we shall omit them.

In 2008, Berg and Nikolaev proved (see \cite[Proposition 14]{BerNik08}) that in any metric space $(X,d)$, the function $\langle\cdot,\cdot\rangle : X^2 \times X^2 \to \mathbb{R}$, defined, for any $x$, $y$, $u$, $v \in X$, by
$$\langle \vv{xy}, \vv{uv} \rangle := \frac12(d^2(x,v) + d^2(y,u) - d^2(x,u) -d^2(y,v))$$
(where an ordered pair of points $(a,b) \in X^2$ is denoted by $\vv{ab}$), called the {\it quasi-linearization function}, is the unique one such that, for any $x$, $y$, $u$, $v$, $w \in X$, we have that:
\begin{enumerate}[(i)]
\item $\langle\vv{xy},\vv{xy}\rangle = d^2(x,y)$;
\item $\langle\vv{xy},\vv{uv}\rangle = \langle\vv{uv},\vv{xy}\rangle$;
\item $\langle\vv{yx},\vv{uv}\rangle = -\langle\vv{xy},\vv{uv}\rangle$;
\item $\langle\vv{xy},\vv{uv}\rangle + \langle\vv{xy},\vv{vw}\rangle = \langle\vv{xy},\vv{uw}\rangle$.
\end{enumerate}
The inner product notation is justified by the fact that if $X$ is a (real) Hilbert space, for any $x$, $y$, $u$, $v \in X$,
\begin{equation}\label{eq-quasi-inner}
\langle\vv{xy},\vv{uv}\rangle = \langle x-y,u-v \rangle = \langle y-x,v-u \rangle.
\end{equation} 
The main result of \cite{BerNik08}, Theorem 1, characterized CAT(0) spaces as being exactly those geodesic spaces $(X,d)$ such that the corresponding Cauchy-Schwarz inequality is satisfied, i.e. for any $x$, $y$, $u$, $v \in X$,
\begin{equation}\label{CauchySchwartz}
\langle\vv{xy},\vv{uv}\rangle \leq d(x,y)d(u,v).
\end{equation}

We shall use, in addition, the following inequality connected to the quasi-linearization function.

\begin{lemma}\label{ineq}
Let $(X,d)$ be a CAT(0) space, $x$, $y$, $z \in X$ and $t \in [0,1]$. Then
$$d^2((1-t)x+ty,z) \leq (1-t)^2d^2(x,z)+2t \langle\vv{yz},\vv{[(1-t)x+ty]z}\rangle.$$
\end{lemma}

\begin{proof}
Set $u:=(1-t)x+ty$. From the defining property of CAT(0) spaces, we have that
$$d^2(z,u) \leq (1-t)d^2(z,x) +td^2(z,y)-t(1-t)d^2(x,y).$$
Multiplying the above by $(1-t)$, and keeping in mind that $d(y,u)=(1-t)d(x,y)$, we have that
\begin{align*}
(1-t)d^2(z,u) &\leq (1-t)^2d^2(z,x) +t(1-t)d^2(z,y)-t(1-t)^2d^2(x,y)\\
&\leq (1-t)^2d^2(z,x) +td^2(z,y)-t(1-t)^2d^2(x,y)\\
&\leq (1-t)^2d^2(z,x) +td^2(z,y)-td^2(y,u),
\end{align*}
so, by adding $td^2(z,u)$, we get that
$$d^2(z,u) \leq (1-t)^2d^2(x,z) + t(d^2(z,y) + d^2(z,u) - d^2(y,u)) = (1-t)^2d^2(x,z) +2t \langle\vv{yz},\vv{uz}\rangle,$$
which is what we needed to show.
\end{proof}

We shall fix now a complete CAT(0) space $(X,d)$ for the remainder of this paper, and throughout the paper, for any self-mapping $T$ of $X$, we shall denote the set of its fixed points by $\Fix(T)$.

A self-mapping $T$ of $X$ is called {\it nonexpansive} if for all $x$, $y \in X$, $d(Tx,Ty) \leq d(x,y)$. If $C \subseteq X$ is closed, convex and nonempty, then there exists a corresponding nearest point projection operator, which one usually denotes by $P_C:X \to C$.

A fundamental property of nonexpansive mappings is the following so-called `resolvent convergence' result, and the original idea of its proof essentially goes back to Minty \cite{Min63}, and was later popularized by Halpern \cite{Hal67}. The generalization to CAT(0) spaces stated below is due to Saejung \cite{Sae10}.

\begin{theorem}[{cf. \cite[Lemmas 2.1 and 2.2]{Sae10}}]\label{s-res}
Let $T:X \to X$ be nonexpansive with $\Fix(T)\neq\emptyset$ and $u \in X$. We have that, for any $t \in (0,1)$ there is a unique $z \in X$ having the property $z=tu+(1-t)Tz$, and we denote it by $z_t$. Then, we have that $\lim_{t \to 0} z_t = P_{\Fix(T)}u$.
\end{theorem}

{\em Firmly nonexpansive} mappings were first introduced, as a refinement of nonexpansive mappings, by Browder \cite{Bro67} in the context of Hilbert spaces and then by Bruck \cite{Bru73} in the context of Banach spaces (this later definition was also studied, e.g., in \cite{Rei77}). The following generalization to geodesic spaces, inspired by the study of firmly nonexpansive mappings in the Hilbert ball \cite{GoeRei82, GoeRei84, ReiSha87, ReiSha90}, was introduced in \cite{AriLeuLop14}.

\begin{definition}
A mapping $T : X \to X$ is called {\em firmly nonexpansive} if for any $x$, $y \in X$ and any $t \in [0,1]$ we have that
$$d(Tx,Ty) \leq d((1-t)x + tTx,(1-t)y+tTy).$$
\end{definition}

As mentioned in \cite{AriLopNic15} (see also \cite{KohLopNic17}), every firmly nonexpansive mapping $T:X \to X$ satisfies the so-called {\em property $(P_2)$}, i.e. that for all $x$, $y \in X$,
$$2d^2(Tx,Ty) \leq d^2(x,Ty) + d^2(y,Tx) - d^2(x,Tx) - d^2(y,Ty),$$
or, using the quasi-linearization function,
\begin{equation}\label{eq-quasi-P2}
d^2(Tx,Ty) \leq \langle\vv{TxTy},\vv{xy}\rangle.
\end{equation}
If $X$ is a Hilbert space, property $(P_2)$ coincides with firm nonexpansiveness as \eqref{eq-quasi-P2} and \eqref{eq-quasi-inner} yield $\|Tx - Ty\|^2 \leq \langle Tx-Ty,x-y \rangle$, which is equivalent to it e.g. by \cite[Proposition 4.2]{BauCom17}. Moreover, from this formulation given by \eqref{eq-quasi-P2} one immediately obtains, using \eqref{CauchySchwartz}, that a self-mapping of a CAT(0) space satisfying property $(P_2)$ is nonexpansive.

Following \cite{LeuNicSip18,Sip21}, if $T$ and $U$ are self-mappings of $X$ and $\lambda$, $\mu>0$, we say that $T$ and $U$ are {\it $(\lambda,\mu)$-mutually firmly nonexpansive} if for all $x$, $y\in X$ and all $\alpha$, $\beta \in [0,1]$ such that 
$(1-\alpha)\lambda=(1-\beta)\mu$, one has that
$$ d(Tx,Uy) \leq d((1-\alpha)x+\alpha Tx,(1-\beta)y+\beta Uy).$$
If $(T_n)_{n \in \N}$ is a family of self-mappings of $X$ and $(\gamma_n)_{n\in\N} \subseteq (0,\infty)$, we say that $(T_n)$ is {\it jointly firmly nonexpansive} with respect to $(\gamma_n)$ if for all $n$, $m\in\N$, $T_n$ and $T_m$ are $(\gamma_n,\gamma_m)$-mutually firmly nonexpansive. In addition, if $(T_\gamma)_{\gamma>0}$ is a family of self-mappings of $X$, we say that it is plainly {\it jointly firmly nonexpansive} if for all $\lambda$, $\mu > 0$, $T_\lambda$ and $T_\mu$ are $(\lambda,\mu)$-mutually firmly nonexpansive. It is clear that a family $(T_\gamma)$ is jointly firmly nonexpansive if and only if for every $(\gamma_n)_{n\in\N} \subseteq (0,\infty)$, $(T_{\gamma_n})_{n \in \N}$ is jointly firmly nonexpansive with respect to $(\gamma_n)$. In \cite{LeuNicSip18} it was shown that examples of jointly firmly nonexpansive families of mappings are furnished by resolvent-type mappings used in convex optimization -- specifically, by:
\begin{itemize}
\item the family $(J_{\gamma f})_{\gamma>0}$, where $f$ is a proper convex lower semicontinous function on $X$ and one denotes for any such function $g$ its proximal mapping by $J_g$;
\item the family $(R_{T, \gamma})_{\gamma>0}$, where $T$ is a nonexpansive self-mapping of $X$ and one denotes, for any $\gamma > 0$, its resolvent of order $\gamma$ by $R_{T, \gamma}$;
\item (if $X$ is a Hilbert space) the family $(J_{\gamma A})_{\gamma>0}$, where $A$ is a maximally monotone operator on $X$ and one denotes for any such operator $B$ its resolvent by $J_B$.
\end{itemize}

Again, if $T$ and $U$ are self-mappings of $X$ and $\lambda$, $\mu>0$, one says that $T$ and $U$ are {\it $(\lambda,\mu)$-mutually $(P_2)$} if for all $x$, $y\in X$, 
$$\frac1{\mu}(d^2(Tx,Uy) + d^2(y,Uy) - d^2(y,Tx)) \leq \frac1{\lambda}(d^2(x,Uy) - d^2(x,Tx) - d^2(Tx,Uy)),$$
or, using the quasi-linearization function,
$$\frac1{\mu}\langle\vv{TxUy},\vv{yUy}\rangle\leq \frac1{\lambda}\langle\vv{TxUy},\vv{xTx}\rangle.$$

\begin{proposition}[{cf. \cite[Proposition 3.10]{LeuNicSip18}}]
Let $\lambda$, $\mu>0$ and $T$ and $U$ be $(\lambda,\mu)$-mutually $(P_2)$ self-mappings of $X$. Then, for all $x \in X$,
$$d(Tx,Ux)\leq \frac{|\lambda-\mu|}\lambda d(x,Tx).$$
\end{proposition}

\begin{corollary}\label{m-n}
Let $\lambda$, $\mu>0$ and $T$ and $U$ be $(\lambda,\mu)$-mutually $(P_2)$ self-mappings of $X$. Then, for all $x \in X$,
$$d(x,Ux)\leq \left(2+\frac\mu\lambda\right) d(x,Tx).$$
\end{corollary}

\begin{proof}
Let $x \in X$. Then
$$d(x,Ux)\leq d(x,Tx) + d(Tx,Ux) \leq d(x,Tx) + \frac{|\lambda-\mu|}\lambda d(x,Tx) \leq  \left(2+\frac\mu\lambda\right) d(x,Tx).$$
\end{proof}

\begin{corollary}[{cf. \cite[Corollary 3.11]{LeuNicSip18}}]\label{same-fixed}
Any two mutually $(P_2)$ self-mappings of $X$ have the same fixed points.
\end{corollary}

One may then similarly state the corresponding definitions for jointly $(P_2)$ families of mappings. As shown in \cite{LeuNicSip18}, all those $(P_2)$ notions generalize their firmly nonexpansive counterparts and coincide with them in the case where $X$ is a Hilbert space. The main result of that paper showed that this condition suffices for the working of the proximal point algorithm, namely that if $X$ is complete, $(T_n)_{n \in \N}$ is a family of self-mappings of $X$ with a common fixed point and $(\gamma_n)_{n\in\N} \subseteq (0,\infty)$ with $\sum_{n=0}^\infty \gamma_n^2 = \infty$, then, assuming that $(T_n)$ is jointly $(P_2)$ with respect to $(\gamma_n)$, any sequence $(x_n) \subseteq X$ such that for all $n$, $x_{n+1}=T_nx_n$, $\Delta$-converges (a generalization of weak convergence to arbitrary metric spaces, due to Lim \cite{Lim76}) to a common fixed point of the family. Moreover, in \cite{Sip21}, the reason for the effectiveness of this sort of condition was further elucidated: Theorem 3.3 of that paper shows that a family of self-mappings is jointly firmly nonexpansive if and only if each mapping in it is nonexpansive and the family as a whole satisfies the well-known resolvent identity.

We will need some facts about sequences of reals. A function $\tau: \N \to \N$ is said to be {\it unboundedly increasing} if $\lim_{n\to\infty}\tau(n)=\infty$ and for all $n \in \N$, $\tau(n) \leq \tau(n+1)$. The following result is immediate.

\begin{lemma}[{\cite[Lemma 2.6]{AoyToy17}}]\label{conv}
Let $(a_n) \subseteq \R$ converging to $0$ and $\tau : \N \to \N$ be unboundedly increasing. Then $\lim_{n \to \infty} a_{\tau(n)} = 0$.
\end{lemma}

\begin{lemma}[{cf. \cite[Lemma 3.1]{Mai08}}]\label{mainge}
Let $(a_n) \subseteq \R$ and $(n_j)$ be a strictly increasing sequence of natural numbers. Assume that for all $j \in \N$, $a_{n_j} < a_{n_j+1}$. Define $\tau: \N \to \N$ by setting, for all $n \in \N$,
$$\tau(n):=\max\{k \leq \max(n_0,n) \mid a_k < a_{k+1}\}.$$
Then:
\begin{itemize}
\item $\tau$ is unboundedly increasing;
\item for all $n \in \N$, $a_{\tau(n)} \leq a_{\tau(n)+1}$ and, for all $n\geq n_0$, $a_n \leq a_{\tau(n)+1}$.
\end{itemize}
\end{lemma}

\begin{corollary}[{\cite[Lemma 2.7]{AoyToy17}}]\label{lat}
Let $(a_n)$ be a non-convergent sequence of nonnegative real numbers. Then there is an $N \in \N$ and an unboundedly increasing $\tau:\N\to\N$ such that for all $n \in\N$, $a_{\tau(n)} \leq a_{\tau(n)+1}$ and, for all $n \geq N$, $a_n \leq a_{\tau(n)+1}$.
\end{corollary}

\begin{proof}
Assume that there is an $n$ such that for all $p>n$, $a_p \geq a_{p+1}$. Then $(a_n)$ is bounded and eventually monotone, so it is convergent, a contradiction. Thus, for all $n$, there is a $p>n$ with $a_p<a_{p+1}$, and by iterating this statement we obtain a sequence $(n_j)$ as in the hypothesis of Lemma~\ref{mainge}. By applying that lemma, we obtain the desired conclusion.
\end{proof}

The following lemma is widely used in fixed point theory.

\begin{lemma}[{\cite[Lemma 2.8]{AoyToy17}}]\label{lfp}
Let $(a_n) \subseteq [0,\infty)$, $(\beta_n) \subseteq \R$ and $(\alpha_n) \subseteq [0,1]$. Suppose that $\sum_{n=0}^\infty \alpha_n = \infty$, $\limsup_{n \to \infty} \beta_n \leq 0$,  and, for all $n$,
$$a_{n+1} \leq (1-\alpha_n)a_n + \alpha_n\beta_n.$$
Then $\lim_{n\to\infty} a_n=0$.
\end{lemma}

We shall also use, for any $a$, $b \in \N$, the notation $[a,b]:=\{n \in \N \mid a \leq n \leq b\}$, disambiguating it by the context from the real interval $[a,b]$ -- we also note that if $a > b$, then $[a,b]=\emptyset$ -- thus one has that $a$, $b\in[a,b]$ only if $a \leq b$, a fact which one must remember to check whenever this property is used.

\section{Convergence theorems}\label{sec:main}

\subsection{Preparatory lemmas}

In this subsection, we state a number of lemmas and propositions which will help us in proving the main convergence theorems, which we do in the next subsection.

The following lemma shows that $(P_2)$ mappings have the property, defined in \cite[Section 4]{Koh16}, of uniform strong quasi-nonexpansiveness, and gives the corresponding `SQNE-modulus'.

\begin{lemma}\label{sne1}
Let $\eps$, $b>0$, $z \in X$, $T:X \to X$ a $(P_2)$ mapping and $p \in \Fix(T)$. Assume that $d(z,p)\leq b$. Then, if
$$d(z,p)-d(Tz,p)<\frac{\eps^2}{2b},$$
we have that $d(z,Tz)<\eps$.
\end{lemma}

\begin{proof}
If $d(z,Tz)=0$, then $d(z,Tz)<\eps$. Assume, then, that $d(z,Tz)\neq0$, so $d(z,p)+d(Tz,p)>0$. Since $T$ is $(P_2)$ and $p \in \Fix(T)$, we have that
$$2d^2(Tz,p) \leq d^2(z,p) + d^2(Tz,p) - d^2(z,Tz),$$
so
$$d^2(Tz,p)\leq d^2(z,p)-d^2(z,Tz).$$
Thus (using, for the strict inequality, the fact that $d(z,p)+d(Tz,p)>0$),
\begin{align*}
d^2(z,Tz) &\leq d^2(z,p) - d^2(Tz,p)\\
&\leq (d(z,p) - d(Tz,p))(d(z,p) + d(Tz,p)) \\
&< \frac{\eps^2}{2b} \cdot (d(z,p) + d(Tz,p)) \\
&\leq \frac{\eps^2}{2b} \cdot 2b = \eps^2,
\end{align*}
so $d(z,Tz) <\eps$.
\end{proof}

\begin{corollary}\label{sne2}
Let $b>0$, $(z_n)\subseteq X$, $(S_n)$ be a family of $(P_2)$ self-mappings of $X$ and $p$ a common fixed point of the $S_n$'s. Assume that, for all $n \in \N$, $d(z_n,p)\leq b$. Then, if
$$\lim_{n \to \infty} (d(z_n,p)-d(S_nz_n,p))=0,$$
we have that $\lim_{n \to \infty} d(z_n,S_nz_n)=0$.
\end{corollary}

\begin{proof}
Let $\eps >0$. We have that there is an $N \in \N$ such that for all $n \geq N$, $d(z_n,p)-d(S_nz_n,p) <\eps^2/(2b)$. Then, by Lemma~\ref{sne1}, for all $n \geq N$, $d(z_n,S_nz_n) <\eps$, from which we get the conclusion.
\end{proof}

\subsection{Main results}

The following lemma, the analogue of \cite[Lemma 2.9]{AoyToy17}, morally forms an integral part of the main convergence proof, so we have chosen to present it in this subsection.

\begin{lemma}\label{l5}
Let $T: X\to X$ be nonexpansive with $\Fix(T)\neq\emptyset$, $(x_n) \subseteq X$ a bounded sequence, $u \in X$ and for all $t \in (0,1)$, let $z_t$ be the unique point in $X$ such that $z_t=tu+(1-t)Tz_t$. Then:
\begin{enumerate}[(i)]
\item \label{p1} for all $t\in(0,1)$ and $n\in\N$, we have that
$$\langle\vv{z_tx_n},\vv{z_tu}\rangle \leq \frac t 2 d^2(x_n,z_t) + \frac{(1-t)^2}{2t}d(x_n,Tx_n)(d(x_n,Tx_n) + 2d(x_n,z_t));$$
\item setting $w:=P_{\Fix(T)}u$, so that $\lim_{t\to 0}z_t = w$, and assuming that $\lim_{n\to\infty} d(x_n,Tx_n)=0$, we have that
$$\limsup_{n\to\infty} \langle\vv{uw},\vv{x_nw}\rangle \leq 0.$$
\end{enumerate}
\end{lemma}

\begin{proof}
\begin{enumerate}[(i)]
\item Using Lemma~\ref{ineq}, we have that
\begin{align*}
d^2(z_t,x_n) &\leq (1-t)^2d^2(Tz_t,x_n) + 2t\langle\vv{ux_n},\vv{z_tx_n}\rangle\\
&\leq (1-t)^2(d(x_n,Tx_n) + d(Tx_n, Tz_t))^2 + 2t(\langle\vv{x_nz_t},\vv{x_nz_t}\rangle + \langle\vv{z_tu},\vv{x_nz_t}\rangle) \\
&\leq (1-t)^2(d^2(x_n,z_t) + d(x_n,Tx_n)(d(x_n,Tx_n)+2d(x_n,z_t))) \\
&\quad + 2t(d^2(x_n,z_t) - \langle\vv{z_tx_n},\vv{z_tu}\rangle),
\end{align*}
from which we get the conclusion.
\item Using \eqref{p1} and that $\lim_{n\to\infty} d(x_n,Tx_n)=0$, we get that for all $t\in(0,1)$.
$$\limsup_{n\to\infty} \langle\vv{z_tx_n},\vv{z_tu}\rangle \leq\frac t 2 \limsup_{n\to\infty} d^2(x_n,z_t).$$
Also, for all $t\in(0,1)$ and $n \in \N$,
\begin{align*}
\prv{uw}{x_nw} &= \prv{uw}{x_nw} - \prv{uw}{x_nz_t} + \prv{uw}{x_nz_t} - \prv{uz_t}{x_nz_t} + \prv{uz_t}{x_nz_t}\\
&= \prv{uw}{x_nw} + \prv{uw}{z_tx_n} + \prv{uw}{x_nz_t} + \prv{z_tu}{x_nz_t} + \prv{z_tx_n}{z_tu}\\
&= \prv{uw}{z_tw} + \prv{z_tw}{x_nz_t} + \prv{z_tx_n}{z_tu}.
\end{align*}
Let $\eps>0$. As $\lim_{t\to 0}z_t = w$, there is a $t_1 \in (0,1)$ such that for all $t\in(0,t_1)$,
$$\prv{uw}{z_tw} \leq d(u,w)d(z_t,w) \leq \frac\eps3,$$
and, using in addition that the set $\{d(x_n,z_t) \mid n \in \N, t \in (0,1)\}$ is bounded (since the curve $(z_t)$ is convergent, hence bounded), we get that there is a $t_2 \in (0,1)$ such that for all $t\in(0,t_2)$ and all $n\in\N$,
$$\prv{z_tw}{x_nz_t}\leq d(z_t,w)d(x_n,z_t) \leq \frac\eps3,$$
and that there is a $t_3 \in (0,1)$ such that for all $t \in (0,t_3)$,
$$\limsup_{n\to\infty} \langle\vv{z_tx_n},\vv{z_tu}\rangle \leq \frac t 2 \limsup_{n\to\infty} d^2(x_n,z_t) \leq \frac\eps3.$$
Let $t\in(0,1)$ be smaller than $t_1$, $t_2$ and $t_3$. Then we get that
$$\limsup_{n\to\infty} \langle\vv{uw},\vv{x_nw}\rangle \leq \eps.$$
As $\eps$ was arbitrarily chosen, we obtain the desired conclusion.
\end{enumerate}
\end{proof}

The following is the main strong convergence theorem of this paper, showing the asymptotic behaviour of the Halpern proximal point algorithm for jointly $(P_2)$ families of mappings.

\begin{theorem}\label{main}
Let $(T_n)$ be a family of self-mappings of $X$, $(\gamma_n) \subseteq (0, \infty)$ and $\gamma>0$ be such that for all $n$, $\gamma_n \geq \gamma$. Assume that the family $(T_n)$ is jointly $(P_2)$ with respect to $(\gamma_n)$. Let $F$ be the common fixed point set of the family and assume that $F\neq\emptyset$. Let $(\alpha_n) \subseteq (0,1]$ such that $\lim_{n\to\infty} \alpha_n = 0$ and $\sum_{n=0}^\infty \alpha_n = \infty$. Let $u \in X$ and $(x_n) \subseteq X$ be such that for all $n$,
$$x_{n+1}=\alpha_n u + (1-\alpha_n)T_nx_n.$$
Then $(x_n)$ converges strongly to $P_Fu$.
\end{theorem}

\begin{proof}
Set $w:=P_Fu$. By Busemann convexity, we have that, for all $n$,
\begin{align*}
d(x_{n+1},w) &\leq \alpha_n d(u,w) + (1-\alpha_n)d(T_nx_n,w) \\
&\leq \alpha_n d(u,w) + (1-\alpha_n)d(x_n,w).
\end{align*}
By induction, one gets that for all $n$, $d(T_nx_n,w)\leq\max(d(u,w),d(x_0,w))$ and thus $(x_n)$ and $(T_nx_n)$ are bounded sequences. Therefore,
$$\lim_{n \to \infty} d(x_{n+1},T_nx_n) = \lim_{n \to \infty} (\alpha_nd(u,T_nx_n)) = 0.$$
Also, we have that, for all $n$,
\begin{align*}
d(x_{n+1},w) &\leq \alpha_n d(u,w) + (1-\alpha_n)d(T_nx_n,w) \\
&\leq \alpha_n d(u,w) + d(T_nx_n,w),
\end{align*}
so, for all $n$,
\begin{equation}\label{min}
d(x_{n+1},w) - d(T_nx_n,w) \leq \alpha_nd(u,w).
\end{equation}
Using Lemma~\ref{ineq} and that, for all $n$, $d(T_nx_n,w)\leq d(x_n,w)$, we have that, for all $n$,
\begin{equation}\label{spec}
d^2(x_{n+1},w) \leq (1-\alpha_n)d^2(x_n,w) + 2\alpha_n\prv{uw}{x_{n+1}w}.
\end{equation}
\noindent {\bf Claim.} The sequence $(d(x_n,w))$ is convergent.\\[1mm]
\noindent {\bf Proof of claim:} Assume towards a contradiction that it is not convergent. Then, by Lemma~\ref{lat}, there is an $N \in \N$ and an unboundedly increasing $\tau:\N\to\N$ such that for all $n \in \N$, $d(x_{\tau(n)},w) \leq d(x_{\tau(n)+1},w)$ and, for all $n \geq N$, $d(x_n,w) \leq d(x_{\tau(n)+1},w)$.

For all $n$, we have that $d(T_{\tau(n)}x_{\tau(n)},w)\leq d(x_{\tau(n)},w)$, so, using \eqref{min}, we get that, for all $n$,
$$0 \leq d(x_{\tau(n)},w)-d(T_{\tau(n)}x_{\tau(n)},w) \leq d(x_{\tau(n)+1},w)-d(T_{\tau(n)}x_{\tau(n)},w)\leq \alpha_{\tau(n)}d(u,w).$$
By Lemma~\ref{conv}, we have that $\lim_{n \to \infty} \alpha_{\tau(n)} = 0$, so, from the above we get that
$$\lim_{n \to \infty} (d(x_{\tau(n)+1},w)-d(T_{\tau(n)}x_{\tau(n)},w)) = 0,$$
and so, by Corollary~\ref{sne2}, that
$$\lim_{n \to \infty} d(T_{\tau(n)}x_{\tau(n)},x_{\tau(n)}) = 0.$$
By Corollary~\ref{m-n}, we have that, for all $n$,
$$d(x_{\tau(n)},T_{\tau(0)}x_{\tau(n)}) \leq \left(2+\frac{\gamma_{\tau(0)}}{\gamma_{\tau(n)}}\right) d(x_{\tau(n)},T_{\tau(n)}x_{\tau(n)}) \leq \left(2+\frac{\gamma_{\tau(0)}}{\gamma}\right) d(x_{\tau(n)},T_{\tau(n)}x_{\tau(n)}),$$
from which we get that
$$\lim_{n \to \infty} d(T_{\tau(0)}x_{\tau(n)},x_{\tau(n)}) = 0.$$
We may now apply Lemma~\ref{l5} to get that
\begin{equation}\label{ls}
\limsup_{n\to\infty} \langle\vv{uw},\vv{x_{\tau(n)}w}\rangle \leq 0.
\end{equation}
On the other hand, we have that, for all $n$,
\begin{align*}
d(x_{\tau(n)},x_{\tau(n)+1}) &\leq d(x_{\tau(n)},T_{\tau(n)}x_{\tau(n)}) + d(T_{\tau(n)}x_{\tau(n)},x_{\tau(n)+1}) \\
&= d(x_{\tau(n)},T_{\tau(n)}x_{\tau(n)}) +\alpha_{\tau(n)}d(u,T_{\tau(n)}x_{\tau(n)}),
\end{align*}
so
$$\lim_{n \to \infty} d(x_{\tau(n)},x_{\tau(n)+1}) = 0.$$
Since, for all $n$,
$$\prv{uw}{x_{\tau(n)}x_{\tau(n)+1}} \leq d(u,w)d(x_{\tau(n)},x_{\tau(n)+1}),$$
we have that
$$\lim_{n \to \infty} \prv{uw}{x_{\tau(n)}x_{\tau(n)+1}} = 0.$$
From the above and \eqref{ls}, we get that
$$\limsup_{n\to\infty} \langle\vv{uw},\vv{x_{\tau(n)+1}w}\rangle \leq 0.$$
Using \eqref{spec}, we have that, for all $n$,
\begin{align*}
d^2(x_{\tau(n)+1},w) &\leq (1-\alpha_{\tau(n)})d^2(x_{\tau(n)},w) + 2\alpha_{\tau(n)}\prv{uw}{x_{\tau(n)+1}w} \\
&\leq (1-\alpha_{\tau(n)})d^2(x_{\tau(n)+1},w) + 2\alpha_{\tau(n)}\prv{uw}{x_{\tau(n)+1}w},
\end{align*}
so, for all $n$,
$$\alpha_{\tau(n)}d^2(x_{\tau(n)+1},w) \leq 2\alpha_{\tau(n)}\prv{uw}{x_{\tau(n)+1}w}.$$
Since, for all $n$, $\alpha_{\tau(n)}>0$, we have that, for all $n$,
$$d^2(x_{\tau(n)+1},w) \leq 2\prv{uw}{x_{\tau(n)+1}w}.$$
Now, for all $n \geq N$, we have that
$$\limsup_{n \to \infty} d^2(x_n,w) \leq \limsup_{n \to \infty}d^2(x_{\tau(n)+1},w) \leq 2\limsup_{n \to \infty}\prv{uw}{x_{\tau(n)+1}w} \leq 0,$$
so $\lim_{n \to \infty} d^2(x_n,w) = 0$, which contradicts our assumption that the sequence $(d(x_n,w))$ is not convergent. This finishes the proof of the claim.\hfill $\blacksquare$ \medbreak

Now, since, for all $n$, $d(T_nx_n,w)\leq d(x_n,w)$, we have that, using \eqref{min},
$$0 \leq d(x_n,w) - d(T_nx_n,w) \leq d(x_n,w) + \alpha_n d(u,w) - d(x_{n+1},w),$$
so
$$\lim_{n \to \infty} (d(x_n,w) - d(T_nx_n,w))=0,$$
and then, by Corollary~\ref{sne2}, that
$$\lim_{n \to \infty} d(x_n,T_nx_n) = 0.$$
By Corollary~\ref{m-n}, we have that, for all $n$,
$$d(x_n,T_0x_n) \leq \left(2+\frac{\gamma_0}{\gamma_n}\right) d(x_n,T_nx_n) \leq \left(2+\frac{\gamma_{0}}{\gamma}\right) d(x_n,T_nx_n),$$
from which we get that
$$\lim_{n \to \infty} d(x_n,T_0x_n)= 0.$$
We may now apply Lemma~\ref{l5} to get that
$$\limsup_{n\to\infty} \langle\vv{uw},\vv{x_nw}\rangle \leq 0,$$
so we also have that
$$\limsup_{n\to\infty} 2\langle\vv{uw},\vv{x_{n+1}w}\rangle \leq 0.$$
By the above, Lemma~\ref{lfp}, and \eqref{spec}, we get that $\lim_{n\to\infty}d^2(x_n,w)=0$ and hence that $\lim_{n\to\infty} x_n=w$.
\end{proof}

The following is the analogue in our context of \cite[Corollary 3.3]{AoyToy17}, giving a convergence theorem for the so-called `Tikhonov regularization' of the proximal point algorithm as discussed in \cite{Xu06}. The fact that convergence results for this kind of iteration may be immediately obtained from the Halpern ones was previously remarked in \cite{LeuNic17}.

\begin{corollary}\label{main-t}
Let $(T_n)$ be a family of self-mappings of $X$, $(\gamma_n) \subseteq (0, \infty)$ and $\gamma>0$ be such that for all $n$, $\gamma_n \geq \gamma$. Assume that the family $(T_n)$ is jointly $(P_2)$ with respect to $(\gamma_n)$. Let $F$ be the common fixed point set of the family and assume that $F\neq\emptyset$. Let $(\beta_n) \subseteq (0,1]$ such that $\lim_{n\to\infty} \beta_n = 0$ and $\sum_{n=0}^\infty \beta_n = \infty$. Let $u \in X$ and $(y_n) \subseteq X$ be such that for all $n$,
$$y_{n+1}=T_n(\beta_n u + (1-\beta_n)y_n).$$
Then $(y_n)$ converges strongly to $P_Fu$.
\end{corollary}

\begin{proof}
For all $n$, put $x_n:=\beta_n u + (1-\beta_n) y_n$ and $\alpha_n:=\beta_{n+1}$. We see that, for all $n$, $y_{n+1}=T_nx_n$ and
$$x_{n+1}=\beta_{n+1}u + (1-\beta_{n+1})y_{n+1} = \alpha_n u + (1-\alpha_n)T_nx_n.$$
We may now apply Theorem~\ref{main} to get that $(x_n)$ converges strongly to $P_Fu$. We also have that, for any $n$,
$$d(y_{n+1},P_Fu) = d(T_nx_n,P_Fu) \leq d(x_n,P_Fu),$$
so $(y_n)$ also converges strongly to $P_Fu$.
\end{proof}

\subsection{Concrete numerical examples}

\begin{figure}
\centering
  \includegraphics[width=220px]{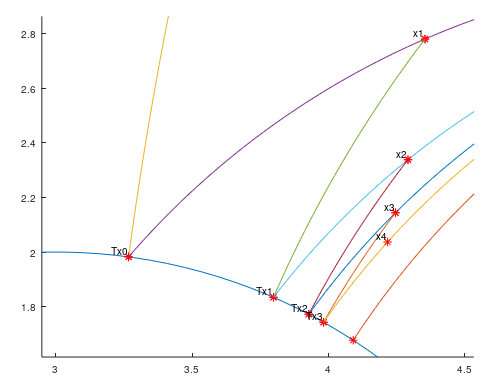}
  \includegraphics[width=220px]{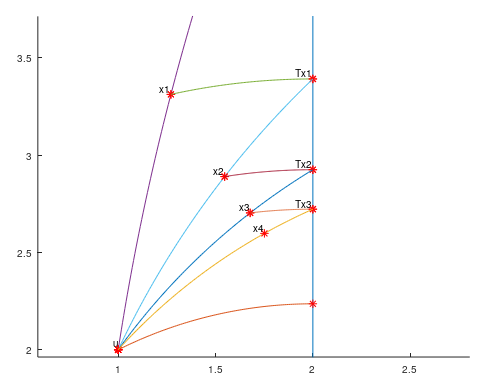}
  \caption{A representation of the first four steps of the iterations.}\label{fig-1}
\end{figure}

As in \cite{AriLopNic15}, we illustrate our main result, Theorem~\ref{main}, with some numerical examples in the upper half-plane model. We wrote an {\it Octave} script to simulate the iteration for the following two sets of inputs: we took the $T_n$'s to be projections onto closed convex sets, namely $\mathcal{C}_{3,2}$ and $\mathcal{R}_2$, respectively, $(\alpha_n)$ to be $(1/(n+2))$ in both cases, $x_0$ to be $(4,5)$ in both cases, and $u$ to be $(6,3)$ and $(1,2)$, respectively. The first few steps of these two simulations can be seen in Figure~\ref{fig-1} and Figure~\ref{fig-2}.

\begin{figure}
\centering
\begin{tabular}{ c|c|c } 
Step & Coordinates & Distance to limit \\
\hline
$x_1$ & (4.354121,   2.781410)&   0.520238\\
$x_2$ & (4.291735,   2.338587)&   0.347794\\
$x_3$ & (4.245949,   2.144022)&   0.259202\\
$x_4$ & (4.216943,   2.036748)&   0.206300\\
$x_5$ & (4.197234,   1.969076)&   0.171248\\
$x_6$ & (4.182998,   1.922578)&   0.146342\\
$x_7$ & (4.172229,   1.888697)&   0.127743\\
\end{tabular}
\begin{tabular}{ c|c|c } 
Step & Coordinates & Distance to limit \\
\hline
$x_1$ & (1.270813,   3.311767)&   0.473128\\
$x_2$ & (1.546908,   2.889523)&   0.311565\\
$x_3$ & (1.677895,   2.702947)&   0.230210\\
$x_4$ & (1.751232,   2.597757)&   0.181888\\
$x_5$ & (1.797667,   2.530419)&   0.150019\\
$x_6$ & (1.829597,   2.483720)&   0.127473\\
$x_7$ & (1.852865,   2.449490)&   0.110707\\
\end{tabular}
    \caption{The coordinates and the distance to the limit point for the first seven iteration steps.}\label{fig-2}
\end{figure}

\section{Quantitative results}\label{sec:quant}

\subsection{Preparatory lemmas}

Similarly to the last section, we present the preparatory lemmas and propositions in a separate subsection.

As stated in the Introduction, if $(x_n)_{n \in \N}$ is a sequence in $X$, then $(x_n)$ is called {\it metastable} if for any $\eps>0$ and $g:\N\to\N$ there is an $N$ such that for all $i$, $j \in [N,N+g(N)]$, $d(x_i,x_j)\leq\eps$, and that a {\it rate of metastability} for $(x_n)$ is a function $\Psi: (0,\infty) \times \N^\N \to \N$ such that for any $\eps$ and $g$, $\Psi(\eps,g)$ gives an upper bound on the (smallest) corresponding $N$. It is an immediate exercise that this is just a reformulation of the Cauchy property.

For all $g: \N \to \N$, we define $\wt{g} : \N \to \N$, for all $n$, by $\wt{g}(n):=n+g(n)$. Also, for all $f:\N \to \N$ and all $n \in \N$, we denote by $f^{(n)}$ the $n$-fold composition of $f$ with itself. Note that for all $g$ and $n$, $\wt{g}^{(n)}(0)\leq\wt{g}^{(n+1)}(0)$. We define, in addition, for any $f :\N \to \N$ and $c \in \N$ the function $f_c: \N \to \N$, setting, for any $l \in \N$, $f_c(l):=f(l+c)$.

The following proposition, which we state in the form that we shall need later, gives a uniform and computable rate of metastability for nonincreasing sequences of nonnegative reals bounded above by a fixed constant.

\begin{proposition}[Quantitative Monotone Convergence Principle, cf. {\cite{Tao08}}]\label{qmcp}
Let $b>0$ and $(a_n)$ be a nonincreasing sequence in $[0,b]$. Then for all $\eps>0$, $g:\N\to\N$ and $l \in \N$ there is an $N \in\left[l, \wt{g}^{\left(\left\lceil\frac b\eps\right\rceil \right)}(l)\right]$ such that for all $i$, $j \in [N,N+g(N)]$, $|a_i-a_j|\leq\eps$.
\end{proposition}

\begin{proof}
Let $\eps>0$ and $g:\N \to \N$. Assume that the conclusion is false, hence in particular for all $i \leq \left\lceil\frac b\eps\right\rceil$, $a_{\wt{g}^{(i)}(l)} - a_{\wt{g}^{(i+1)}(l)} > \eps$. Then
$$b \geq a_l \geq a_l -  a_{\wt{g}^{\left(\left\lceil\frac b\eps\right\rceil + 1\right)}(l)} = \sum_{i=0}^{\left\lceil\frac b\eps\right\rceil} \left(a_{\wt{g}^{(i)}(l)} - a_{\wt{g}^{(i+1)}(l)}\right) > \left\lceil\frac b\eps\right\rceil\cdot\eps \geq b,$$
a contradiction.
\end{proof}

The following argument is implicit in the proof of \cite[Theorem 4.1]{Koh20}.

\begin{lemma}\label{l2b}
Let $b>0$ and $x$, $y$, $z \in X$ be such that $d(x,y) \leq b$ and $d(y,z)\leq b$. Then
$$d^2(x,y)\leq d^2(y,z) + 2bd(x,z).$$
\end{lemma}

\begin{proof}
Since $d(x,y) \leq d(x,z) + d(y,z)$, we have that $d(x,y)-d(y,z)\leq d(x,z)$. Now,
$$d^2(x,y) - d^2(y,z) = (d(x,y)+d(y,z))(d(x,y) - d(y,z)) \leq (d(x,y)+d(y,z))\cdot d(x,z) \leq 2b \cdot d(x,z),$$
from which the conclusion follows.
\end{proof}

We shall now present the analogue in our context of \cite[Lemma 3.3]{Koh20}, which is proven there using full strong nonexpansiveness. Since the adaptation of that concept to the metric context would be `somewhat artificial' \cite[p. 229]{Koh16}, we are being led, as we already said in the Introduction, to further mine the already partially quantitative Lemma~\ref{sne1} into the following property which we then dubbed `quantitative quasiness'.

\begin{proposition}\label{sne-q}
Denote, for this and subsequent results, for any $\eps$, $b>0$, $\omega(b,\eps):=\frac{\eps^2}{15b}$. Let $\eps$, $b>0$, $z$, $p\in X$ and $T: X\to X$ a $(P_2)$ mapping. Assume that $d(z,p) \leq b$ and $d(p,Tp) \leq b$. Then, if
$$d(z,p)-d(Tz,p) \leq \omega(b,\eps)$$
and
$$d(p,Tp)\leq \omega(b,\eps),$$
we have that $d(z,Tz)\leq\eps$.
\end{proposition}

\begin{proof}
Since $T$ is $(P_2)$, we have that
\begin{align*}
2d^2(Tz,Tp) &\leq d^2(z,Tp) + d^2(Tz,p) - d^2(z,Tz) - d^2(p,Tp)\\
& \leq (d(z,p)+d(p,Tp))^2 + d^2(Tz,p) - d^2(z,Tz) - d^2(p,Tp).
\end{align*}
As
$$d(Tz,Tp) \geq |d(Tz,p) - d(p,Tp)|,$$
we have that
$$d^2(Tz,Tp) \geq d^2(Tz,p) + d^2(p,Tp) - 2d(Tz,p)d(p,Tp) \geq d^2(Tz,p)  - 2d(Tz,p)d(p,Tp),$$
so
$$2d^2(Tz,p)  - 4d(Tz,p)d(p,Tp) \leq d^2(z,p) + d^2(p,Tp) + 2d(z,p)d(p,Tp) + d^2(Tz,p) - d^2(z,Tz) - d^2(p,Tp).$$
Thus,
\begin{align*}
d^2(z,Tz) &\leq d^2(z,p) - d^2(Tz,p) + 4d(p,Tp)(d(z,p)+d(Tz,p))\\
& = (d(z,p)-d(Tz,p)+4d(p,Tp))(d(z,p)+d(Tz,p))\\
& \leq \left(\frac{\eps^2}{15b} + 4 \cdot \frac{\eps^2}{15b}\right)(d(z,p)+d(Tz,p))\\
& = \frac{\eps^2}{3b} \cdot(d(z,p)+d(Tz,p))\\
& \leq \frac{\eps^2}{3b} \cdot(d(z,p)+d(Tz,Tp)+d(p,Tp))\\
& \leq \frac{\eps^2}{3b} \cdot(2d(z,p)+d(p,Tp))\\
& \leq \frac{\eps^2}{3b} \cdot(2b+b) = \frac{\eps^2}{3b} \cdot 3b = \eps^2,
\end{align*}
so $d(z,Tz) \leq\eps$.
\end{proof}

The following is the quantitative version of Lemma~\ref{l5}, i.e. the analogue of \cite[Lemma 3.6]{Koh20}.

\begin{lemma}\label{l5q}
Let $T: X\to X$ be nonexpansive with $\Fix(T)\neq\emptyset$, $(x_n) \subseteq X$ a bounded sequence, $u \in X$ and for all $t \in (0,1)$, let $z_t$ be the unique point in $X$ such that $z_t=tu+(1-t)Tz_t$. Let $b>0$ such that for all $n \in \N$ and all $t\in (0,1)$ one has $d(z_t,x_n) \leq b$ and $d(x_n,Tx_n) \leq b$. Let $(t_l)_{l \in \N^*} \subseteq (0,1)$ and $\rho:(0,\infty)\to\N$ be such that for all $l\geq \rho(\eps)$ we have that $t_l \leq\eps$. Let $\chi : \N^* \to \N^*$ be such that for all $l \in \N^*$ we have that $t_l \geq \frac1{\chi(l)}$, i.e. $\frac1{t_l\chi(l)} \leq 1$. Take $k \geq \rho\left(\frac\eps{b^2}\right)$, so that $\frac{t_k}2 \cdot b^2 \leq \frac\eps2$. Take $n$ such that $d(x_n,Tx_n) \leq \frac{\eps}{3b\chi(k)}$. Then
$$\prv{uz_{t_k}}{x_nz_{t_k}} \leq \eps.$$
\end{lemma}

\begin{proof}
By Lemma~\ref{l5}.\eqref{p1}, we have that
\begin{align*}
\prv{uz_{t_k}}{x_nz_{t_k}} &\leq \frac {t_k} 2 d^2(x_n,z_{t_k}) + \frac{(1-{t_k})^2}{2{t_k}}d(x_n,Tx_n)(d(x_n,Tx_n) + 2d(x_n,z_{t_k}))\\
&\leq \frac{t_k}2 \cdot b^2 + \frac{3b}{2t_k}d(x_n,Tx_n)\\
&\leq \frac\eps2 + \frac{3b}{2t_k}\cdot\frac{\eps}{3b\chi(k)} \leq \frac\eps2+\frac\eps2 = \eps.
\end{align*}
\end{proof}

\begin{lemma}[{\cite[Lemma 3.4.1]{Koh20}}]\label{l91}
For any $\eps>0$, $g: \N \to \N$, $K\in\N$, $b>0$, set
$$\psi(\eps,g,K,b):=\wt{g}^{\left(\left\lceil\frac{b}{\eps}\right\rceil\right)}(K) \geq K.$$
Let $b>0$, $(a_n) \subseteq [0,b]$ and $\tau:\N\to\N$ such that for all $k$, $n\in\N$ with $k \leq n$ and $a_k< a_{k+1}$, we have $k \leq \tau(n)$.

Then, for all $g: \N\to\N$, $K\in\N$ and $\eps>0$ with $\tau(\psi(\eps,g,K,b)) <K$ we have that there is a $n \in [K,\psi(\eps,g,K,b)]$ such that for all $i$, $j \in [n,n+g(n)]$, $|a_i-a_j|\leq\eps$.
\end{lemma}

The following is the quantitative version of Lemma~\ref{lfp}.

\begin{lemma}[{\cite[Lemma 3.5]{Koh20}}]\label{l10}
For any $\eps>0$, $S:(0,\infty) \times \N \to \N$, $m\in\N$ and $b>0$, set
$$\vp(\eps,S,m,b):=m+S\left(\frac{\eps}{4b},m\right)+1.$$
Let $b>0$, $(a_n)\subseteq [0,b]$, $(\alpha_n) \subseteq (0,1]$, $(\beta_n)\subseteq \R$ and $(\gamma_n)\subseteq[0,\infty)$. Suppose that for any $n \in\N$,
$$a_{n+1}\leq (1-\alpha_n)a_n+\alpha_n\beta_n+\gamma_n.$$
Let $S:(0,\infty)\times \N \to \N$ be nondecreasing in the second argument such that for all $\eps>0$ and $m \in \N$,
$$\prod_{k=m}^{S(\eps,m)} (1-\alpha_k) \leq\eps.$$
Let $\eps>0$, $g: \N \to \N$ and $P \in \N$ be such that there is an $m\leq P$ such that for all
$$i \in \left[m,m+g^M\left(m+S\left(\frac\eps{4b},m\right)+1\right)+S\left(\frac\eps{4b},m\right)\right],$$
we have that $\beta_i\leq\frac\eps4$. Suppose that
$$\sum_{i=0}^{\vp(\eps,S,P,b)+g^M(\vp(\eps,S,P,b))} \gamma_i \leq \frac\eps2.$$

Then there is an $N \leq \vp(\eps,S,P,b)$ such that for all $i \in [N,N+g(N)]$, $a_i\leq\eps$.
\end{lemma}

\begin{lemma}[{\cite[Lemma 3.7]{Koh20}}]\label{l12}
Let $(y_n)_{n \geq 1} \subseteq X$ and $\xi:(0,\infty) \times \N^\N \to \N$ be such that for any $\eps>0$ and $g: \N\to\N$ there is an $n\in[1,\xi(\eps,g)]$ such that for all $i$, $j \in [n,g(n)]$, $d(y_i,y_j)\leq\eps$.

Then there is an $\eps>0$ such that for all $c\in\N^*$ and all $f:\N\to\N$ there is a $k \in [c,\xi(\eps,f_c)+c]$ such that for all $i$, $j\in[k,f(k)]$, $d(y_i,y_j)\leq\eps$.
\end{lemma}

The following is the quantitative version of Theorem~\ref{s-res}, as obtained in \cite{KohLeu12}. An abstract version of it which uses the concept of jointly firmly nonexpansive families of mappings may be found in \cite[Section 5]{Sip21}, but here we shall only need the rate of metastability for the resolvents of nonexpansive mappings.

\begin{proposition}[{cf. \cite[Proposition 9.3]{KohLeu12}}]\label{l-res}
Define, for all $b$, $\eps>0$ and $g:\N\to\N$, $\xi_b(\eps,g):=g^{\left(\left\lceil\frac{b^2}{\eps^2}\right\rceil\right)}(1)$.

Let $T:X\to X$ be nonexpansive, $u \in X$, and for all $t \in (0,1)$ put $z_t$ to be the unique point in $X$ such that $z_t=tu+(1-t)Tz_t$. Let $(t_n)_{n \in \N^*} \subseteq [0,1]$ be nonincreasing. Put, for any $n \in \N^*$, $y_n:=z_{t_n}$. Let $b>0$ and assume that, for all $n$, $d(y_n,u)\leq b$. Then, for any $\eps>0$ and any $g:\N\to\N$ there is an $n \leq \xi_b(\eps,g)$ such that for all $i$, $j \in [n,g(n)]$, $d(y_i,y_j)\leq\eps$.
\end{proposition}

\subsection{Main results}

The main quantitative theorem includes, as expected, a rate of metastability, and in order to express it we shall introduce the following notations.

\begin{notation}\label{not}
Let $b$, $\gamma>0$, $(\wt{\gamma}_n) \subseteq (0,\infty)$, $(\wt{\alpha}_n) \subseteq (0,1]$, $\zeta:(0,\infty) \to \N$, $S:(0,\infty)\times \N \to \N$, $\eps>0$ and $g:\N\to\N$.

We shall introduce a series of quantities depending on these parameters. Set
$$C:=2+\frac{\wt{\gamma}_0}\gamma, \quad \hat{\eps}:=\frac{\eps^2}{128b}.$$
Set, for all $l \in \N$,
$$\eta_l:=\frac{\eps^2}{192bl},\quad M_1(l):=\min\left(\frac12\omega\left(b,\frac{\eta_l}C\right),\omega\left(b,\frac{\eps^2}{128b}\right),\frac{\eps^2}{128b}\right),\quad n_l:=\max\left\{\zeta\left(\frac{M_1(i)}{b}\right) \bigg|\ i\leq l \right\},$$
$$\hat{g}(l):=g^M\left(l+S\left(\frac{\eps^2}{16b^2},l\right)+1\right)+S\left(\frac{\eps^2}{16b^2},l\right),\quad g'(l):=\hat{g}(l)+2.$$
Set, for all $l$, $i \in \N$,
$$\theta(l,i):= \psi\left(\frac12\omega\left(b,\frac{\eta_l}C\right),g',i,b\right)\geq i,$$
and for all $l \in \N$,
$$\theta^*(l):=\max\{\theta(j,n_j) \mid j \leq l\},\quad K(l):=\theta(l,n_l)+\hat{g}^M(\theta(l,n_l))+2,$$
$$\hat{K}(l):=K(l)+S\left(\frac{\eps^2}{16b^2},K(l)\right)+1+g^M\left(K(l)+S\left(\frac{\eps^2}{16b^2},K(l)\right)+1\right),$$
$$\wt{\gamma}^M_l:=\max\{\wt{\gamma}_j \mid j\leq l\}.$$
Set, for all $\beta >0$ and $l \in \N$,
$$\wt{\rho}(\beta,l):=\left\lceil \frac{\left(2+\frac{\wt{\gamma}^M_l}\gamma\right) \cdot b}{\beta}\right\rceil.$$
Set, for all $l \in \N$,
$$M_2(l):=\min\left\{\frac\eps2,\frac{\eps^2}{16b\left(\hat{K}(l)+1\right)},\omega\left(b,\frac{\eps^2}{128b}\right),\omega\left(b,\frac{\eta_l}C\right),\frac{\eps^2}{16b}\cdot \min\{\wt{\alpha}_j \mid j \leq K(l) \}\right\},$$
$$f(l):=\max\left(\wt{\rho}\left(M_2(l),\hat{K}(l)\right),l\right)\geq l.$$
Set, now,
$$c:=\left\lceil\frac{64b^2}{\eps^2}\right\rceil,\quad k^*:=\xi_b(\hat{\eps},f_c)+c,\quad K^*:=\theta^*(k^*) +\hat{g}^M(\theta^*(k^*))+2,$$
$$\Phi:=K^*+S\left(\frac{\eps^2}{16b^2},K^*\right)+1.$$
This last quantity we shall denote in the sequel by $\Phi_{b,\gamma,(\wt{\gamma}_n),(\wt{\alpha}_n),\zeta,S}(\eps,g)$, i.e. explicitly expressing its dependence on the parameters.
\end{notation}

Armed with the above, we may now state the quantitative version of Theorem~\ref{main}, which gives a rate of metastability for the Halpern proximal point algorithm in our context.

\begin{theorem}\label{main-q}
Let $(T_n)$ be a family of self-mappings of $X$, $(\gamma_n) \subseteq (0,\infty)$ and $\gamma>0$ be such that $(T_n)$ is jointly $(P_2)$ with respect to $(\gamma_n)$ and for all $n$, $\gamma_n \geq \gamma$. Let $(\wt{\gamma}_n) \subseteq (0,\infty)$ be such that for all $n$, $\wt{\gamma}_n \geq \gamma_n$. We denote by $F$ the common fixed point set of the family $(T_n)$. Let $(\alpha_n) \subseteq (0,1]$, $u\in X$ and $(x_n)\subseteq X$ be such that for all $n$,
$$x_{n+1}=\alpha_n u +(1-\alpha_n)T_nx_n.$$
Let $(\wt{\alpha}_n) \subseteq (0,1]$ be such that for all $n$, $\wt{\alpha}_n \leq \alpha_n$. Let $b \in \N^*$ and $p \in F$ be such that $2d(x_0,p) \leq b$ and $2d(u,p)\leq b$. Let $\zeta:(0,\infty) \to \N$ be such that for all $\beta>0$ and all $m \geq \zeta(\beta)$, $\alpha_m\leq\beta$. Let $S:(0,\infty)\times \N \to \N$ be nondecreasing in the second argument such that for all $\eps>0$ and $m \in \N$,
$$\prod_{k=m}^{S(\eps,m)} (1-\alpha_k) \leq\eps.$$
Then:
\begin{enumerate}[(i)]
\item for all $\eps>0$ and $g: \N \to \N$, there is a $w \in X$ and an $N \leq \Phi_{b,\gamma,(\wt{\gamma}_n),(\wt{\alpha}_n),\zeta,S}(\eps,g)$ such that for all $i \in [N,N+g(N)]$, $d(w,T_iw) \leq \eps/2$ and $d(x_i,w)\leq\eps/2$.
\item $\Phi_{b,\gamma,(\wt{\gamma}_n),(\wt{\alpha}_n),\zeta,S}$ is a rate of metastability for $(x_n)$, i.e. for all $\eps>0$ and $g: \N \to \N$, there is an $N \leq \Phi_{b,\gamma,(\wt{\gamma}_n),(\wt{\alpha}_n),\zeta,S}(\eps,g)$ such that for all $i$, $j \in [N,N+g(N)]$, $d(x_i,x_j)\leq\eps$.
\end{enumerate} 
\end{theorem}

\begin{proof}
Let $\eps>0$ and $g:\N\to\N$. We shall use the notations from Notation~\ref{not}, instantiating the parameters with those from the statement of the theorem, together with this $\eps$ and $g$.

We first remark that the second bullet point is an immediate consequence of the first one.

For all $t\in(0,1)$, set $z_t$ to be the unique point such that $z_t=tu+(1-t)T_0z_t$. Note that, for all $t \in (0,1)$, by Busemann convexity, we have that
$$d(z_t,p) \leq td_t(u,p) + (1-t)d(T_0z_t,p) \leq td_t(u,p) + (1-t)d(z_t,p),$$
so $d(z_t,p) \leq d(u,p)$, from which we get $d(z_t,u) \leq 2d(u,p) \leq b$. Also, for all $n$ and $t$, $d(T_nz_t,p) \leq d(z_t,p) \leq d(u,p)$, so $d(T_nz_t,u) \leq 2d(u,p) \leq b$ and $d(T_nz_t,z_t) \leq d(T_nz_t,p)+d(z_t,p)\leq 2(u,p) \leq b$.

Again by Busemann convexity, we have that, for all $n$,
\begin{align*}
d(x_{n+1},p) &\leq \alpha_n d(u,p) + (1-\alpha_n)d(T_nx_n,p) \\
&\leq \alpha_n d(u,p) + (1-\alpha_n)d(x_n,p).
\end{align*}
By induction, one gets that for all $n$, $d(T_nx_n,p)\leq d(x_n,p)\leq \max(d(u,p),d(x_0,p))\leq b/2$, so, for all $n$, $d(T_nx_n,u) \leq d(T_nx_n,p)+d(u,p)\leq b$ and, for all $n$ and $t$, $d(x_n,z_t)  \leq d(x_n,p) + d(z_t,p) \leq b$.

For all $l \in \N^*$, set $y_l:=z_{1/l}$. By Proposition~\ref{l-res} and Lemma~\ref{l12}, we get that there is a $k \in [c,k^*]$ such that for all $i$, $j \in [k,f(k)]$, $d(y_i,y_j)\leq\hat{\eps}$. Set $k':=f(k)\geq k$. We get in particular that $d(y_k,y_{k'})\leq\hat{\eps}$ and that
$$d(y_{k'},T_0y_{k'}) = d\left(\frac1{k'}u + \left(1-\frac1{k'}\right)T_0y_{k'},T_0y_{k'}\right) = \frac1{k'}d(u,T_0y_{k'}) \leq \frac{b}{k'}.$$

We shall take $w:=y_{k'}$ and thus it remains to be shown that there is an $N \leq \Phi$ such that for all $i \in [N,N+g(N)]$, $d(y_{k'},T_iy_{k'}) \leq\eps/2$ and $d(x_i,y_{k'})\leq\eps/2$.

Set $A:=\theta(k,n_k)\geq k$ and, for all $m$, $a_m:=d(x_m,y_{k'})$. We distinguish two cases.

{\bf Case I.} For all $i \leq A$, $a_{i+1}\leq a_i$.

Since
$$\theta(k,n_k) = \psi\left(\frac12\omega\left(b,\frac{\eta_k}C\right),g',n_k,b\right)=\wt{g'}^{\left(\left\lceil \frac{b}{\frac12\omega\left(b,\frac{\eta_k}C\right)} \right\rceil \right)}(n_k),$$
we get by Proposition~\ref{qmcp} that there is an $n\in [n_k,A]$ such that for all $i$, $j \in [n,n+g'(n)]=[n,n+\hat{g}(n)+2]$, $|a_i-a_j|\leq\frac12\omega\left(b,\frac{\eta_k}C\right)$. We keep this in mind.

{\bf Case II.} There is an $i \leq A$ with $a_{i+1}> a_i$.

Define $\tau:\N\to\N$, for all $n \in \N$, by
$$\tau(n):=\max\{j \leq \max(n,A) \mid a_j < a_{j+1}\}.$$
Then:
\begin{itemize}
\item for all $n \in \N$, $\tau(n)\leq\tau(n+1)$ and $a_{\tau(n)} \leq a_{\tau(n)+1}$;
\item for all $l$, $n\in\N$ with $l \leq n$ and $a_l< a_{l+1}$, we have $l \leq \tau(n)$;
\item for all $n\geq A$, $a_n \leq a_{\tau(n)+1}$ (this is the only non-trivial statement, but \cite[Lemma 3.4.2]{Koh20} shows that it follows exactly as in the original proof of Lemma~\ref{mainge}, i.e. see \cite[Lemma 3.1]{Mai08}).
\end{itemize}

We now distinguish two sub-cases.

{\bf Sub-case II.1.} For all $m \in [A,A+\hat{g}(A)+2]$, $\tau(m)\geq n_k$.

Let $m \in [A,A+\hat{g}(A)+2]$ be arbitrary. Then
$$d(x_{\tau(m)+1},y_{k'}) \leq \alpha_{\tau(m)}d(u,y_{k'})+(1-\alpha_{\tau(m)})d(T_{\tau(m)}x_{\tau(m)},y_{k'}) \leq \alpha_{\tau(m)}d(u,y_{k'})+d(T_{\tau(m)}x_{\tau(m)},y_{k'}),$$
so
$$d(x_{\tau(m)+1},y_{k'})-d(T_{\tau(m)}x_{\tau(m)},y_{k'}) \leq \alpha_{\tau(m)}d(u,y_{k'}) \leq \alpha_{\tau(m)}b$$
and (using that $\tau(m)\geq n_k$)
\begin{align*}
d(x_{\tau(m)},y_{k'})-d(T_{\tau(m)}x_{\tau(m)},y_{k'}) &\leq d(x_{\tau(m)+1},y_{k'})-d(T_{\tau(m)}x_{\tau(m)},y_{k'}) \\
&\leq \alpha_{\tau(m)}b \leq M_1(k) \leq \min\left(\omega\left(b,\frac{\eta_k}C\right),\omega\left(b,\frac{\eps^2}{128b}\right)\right).
\end{align*}

As
$$\tau(m) \leq \max(m,A) =m \leq A+\hat{g}(A)+2\leq A+\hat{g}^M(A)+2 = K(k) \leq \hat{K}(k),$$
we have that
\begin{align*}
d(y_{k'},T_{\tau(m)}y_{k'}) &\leq \left( 2+\frac{\gamma_{\tau(m)}}{\gamma_0} \right)\cdot d(y_{k'},T_0y_{k'}) \leq \left( 2+\frac{\wt{\gamma}^M_{ \hat{K}(k)}}{\gamma} \right)\cdot \frac{b}{k'} \\
&\leq \left( 2+\frac{\wt{\gamma}^M_{ \hat{K}(k)}}{\gamma} \right) \cdot b \cdot \frac1{\left\lceil \frac{\left( 2+\frac{\wt{\gamma}^M_{ \hat{K}(k)}}{\gamma} \right) \cdot b}{M_2(k)}\right\rceil}  \leq M_2(k) \leq \min\left(\omega\left(b,\frac{\eta_k}C\right),\omega\left(b,\frac{\eps^2}{128b}\right)\right).
\end{align*}
By Proposition~\ref{sne-q}, we get that
$$d(x_{\tau(m)},T_{\tau(m)}x_{\tau(m)}) \leq \min\left(\frac{\eta_k}C,\frac{\eps^2}{128b}\right),$$
so
$$d(x_{\tau(m)},T_0x_{\tau(m)}) \leq \left( 2+\frac{\gamma_0}{\gamma_{\tau(m)}} \right) d(x_{\tau(m)},T_{\tau(m)}x_{\tau(m)}) \leq \left( 2+\frac{\wt{\gamma}_0}{\gamma} \right)\cdot \frac{\eta_k}C = \eta_k,$$
and, since $\tau(m) \geq n_k$,
\begin{align*}
d(x_{\tau(m)+1},x_{\tau(m)}) &\leq d(x_{\tau(m)+1},T_{\tau(m)}x_{\tau(m)}) + d(T_{\tau(m)}x_{\tau(m)},x_{\tau(m)})\\
&= \alpha_{\tau(m)} d(u,T_{\tau(m)}x_{\tau(m)}) + d(T_{\tau(m)}x_{\tau(m)},x_{\tau(m)})\\
&\leq \alpha_{\tau(m)}\cdot b+ \frac{\eps^2}{128b} \leq \frac{\eps^2}{128b}+\frac{\eps^2}{128b}=\frac{\eps^2}{64b}.
\end{align*}

As $d(x_{\tau(m)},T_0x_{\tau(m)}) \leq\eta_k$ and $k \geq c = \left\lceil\frac{64b^2}{\eps^2}\right\rceil$, we have, by Lemma~\ref{l5q}, that
$$\prv{uy_k}{x_{\tau(m)}y_k} \leq \frac{\eps^2}{64}.$$
On the other hand,
$$\prv{uy_k}{x_{\tau(m)}x_{\tau(m)+1}} \leq d(u,y_k)d(x_{\tau(m)},x_{\tau(m)+1}) \leq b \cdot \frac{\eps^2}{64b}=\frac{\eps^2}{64},$$
so $\prv{uy_k}{x_{\tau(m)+1}y_k} \leq \frac{\eps^2}{32}$. We also know that $d(y_k,y_{k'}) \leq \hat{\eps} = \frac{\eps^2}{128b}$, so
\begin{align*}
\prv{uy_{k'}}{x_{\tau(m)+1}y_{k'}} &= \prv{uy_k}{x_{\tau(m)+1}y_k} + \prv{uy_k}{y_ky_{k'}} + \prv{y_ky_{k'}}{x_{\tau(m)+1}y_k'}\\
&\leq \frac{\eps^2}{32} + \frac{\eps^2}{128} + \frac{\eps^2}{128} < \frac{\eps^2}{16}.
\end{align*}
Using Lemma~\ref{ineq} and Lemma~\ref{l2b}, we get that
\begin{align*}
d^2(x_{\tau(m)+1},y_k) &\leq (1-\alpha_{\tau(m)})^2d^2(T_{\tau(m)}x_{\tau(m)},y_{k'}) + 2\alpha_{\tau(m)} \prv{uy_{k'}}{x_{\tau(m)+1}y_{k'}}\\
&\leq (1-\alpha_{\tau(m)})^2d^2(T_{\tau(m)}x_{\tau(m)},T_{\tau(m)}y_{k'}) + 2bd(y_{k'},T_{\tau(m)}y_{k'}) + 2\alpha_{\tau(m)} \prv{uy_{k'}}{x_{\tau(m)+1}y_{k'}}\\
&\leq (1-\alpha_{\tau(m)})d^2(x_{\tau(m)},y_{k'}) + 2bd(y_{k'},T_{\tau(m)}y_{k'}) + 2\alpha_{\tau(m)} \prv{uy_{k'}}{x_{\tau(m)+1}y_{k'}}\\
&\leq (1-\alpha_{\tau(m)})d^2(x_{\tau(m)+1},y_{k'}) + 2bd(y_{k'},T_{\tau(m)}y_{k'}) + 2\alpha_{\tau(m)} \prv{uy_{k'}}{x_{\tau(m)+1}y_{k'}},
\end{align*}
so
$$d^2(x_{\tau(m)+1},y_{k'}) \leq 2\prv{uy_{k'}}{x_{\tau(m)+1}y_{k'}} + \frac{2bd(y_{k'},T_{\tau(m)}y_{k'})}{\alpha_{\tau(m)}} \leq \frac{\eps^2}8 +\frac{\eps^2}8 \leq \frac{\eps^2}4.$$

Since $m \geq A$, we have that $d^2(x_m,y_{k'}) \leq d^2(x_{\tau(m)+1},y_{k'}) \leq {\eps^2}/4$, so $d(x_m,y_{k'}) \leq\eps/2$.

As $m \leq A+\hat{g}(A)+2 = K(k) \leq \hat{K}(k)$, we have that
\begin{align*}
d(y_{k'},T_my_{k'}) &\leq \left( 2+\frac{\gamma_{m}}{\gamma_0} \right)\cdot d(y_{k'},T_0y_{k'}) \leq \left( 2+\frac{\wt{\gamma}^M_{ \hat{K}(k)}}{\gamma} \right)\cdot \frac{b}{k'} \\
&\leq \left( 2+\frac{\wt{\gamma}^M_{ \hat{K}(k)}}{\gamma} \right) \cdot b \cdot \frac1{\left\lceil \frac{\left( 2+\frac{\wt{\gamma}^M_{ \hat{K}(k)}}{\gamma} \right) \cdot b}{M_2(k)}\right\rceil}  \leq M_2(k) \leq \frac\eps2.
\end{align*}

As $m$ was arbitrarily chosen, we have shown that for all $m \in [A,A+\hat{g}(A)+2]$, $d(y_{k'},T_my_{k'})\leq\eps/2$ and $d(x_m,y_{k'}) \leq\eps/2$.

We can then take $N:=A$, because then, as $A=\theta(k,n_k)\leq K(k)$ and $k \leq k^*$, we have $\theta(k,n_k) \leq \theta^*(k^*)$ and so $\hat{g}^M(\theta(k,n_k)) \leq \hat{g}^M(\theta^*(k^*))$, $\theta(k,n_k)+\hat{g}^M(\theta(k,n_k))+2 \leq \theta^*(k^*)+\hat{g}^M(\theta^*(k^*))+2$, so $K(k)\leq K^* \leq \Phi$. We have thus shown $N \leq \Phi$ and we derive the needed conclusion by noting that $g(N) \leq \hat{g}(N)$.

{\bf Sub-case II.2.} There is an $m \in [A,A+\hat{g}(A)+2]$ with $\tau(m)< n_k$.

Since $A\leq m$, $\tau(A)\leq\tau(m)<n_k$. But $\theta(k,n_k) = \psi\left(\frac12\omega\left(b,\frac{\eta_k}C\right),g',n_k,b\right)$, so, by Lemma~\ref{l91}, we get that there is an $n\in [n_k,A]$ such that for all $i$, $j \in [n,n+g'(n)]=[n,n+\hat{g}(n)+2]$, $|a_i-a_j|\leq\frac12\omega\left(b,\frac{\eta_k}C\right)$. We note that this was also proven in Case I, so now we may merge the two threads of the proof (and we no longer need this $m$ above).

Note that, since $n \leq A$, we have that $n+g'(n) \leq A+(g')^M(A) = K(k) \leq \hat{K}(k)$. Also note that, as $n \geq n_k$, for all $m \geq n$, $\alpha_m b\leq M_1(k) \leq \frac12\omega\left(b,\frac{\eta_k}C\right)$.

Let $m \in [n,n+g'(n)-1]=[n,n+\hat{g}(n)+1]$. We have that
\begin{align*}
d(x_m,y_{k'}) - d(T_mx_m,y_{k'}) &= d(x_{m+1},y_{k'}) - d(T_mx_m,y_{k'}) +d(x_m,y_{k'}) -d(x_{m+1},y_{k'})\\
&\leq \alpha_m b  + \frac12\omega\left(b,\frac{\eta_k}C\right) \leq\omega\left(b,\frac{\eta_k}C\right).
\end{align*}
As $m \leq n+g'(n) \leq \hat{K}(k)$, we have that
\begin{align*}
d(y_{k'},T_my_{k'}) &\leq \left( 2+\frac{\gamma_{m}}{\gamma_0} \right)\cdot d(y_{k'},T_0y_{k'}) \leq \left( 2+\frac{\wt{\gamma}^M_{ \hat{K}(k)}}{\gamma} \right)\cdot \frac{b}{k'} \\
&\leq \left( 2+\frac{\wt{\gamma}^M_{ \hat{K}(k)}}{\gamma} \right) \cdot b \cdot \frac1{\left\lceil \frac{\left( 2+\frac{\wt{\gamma}^M_{ \hat{K}(k)}}{\gamma} \right) \cdot b}{M_2(k)}\right\rceil}  \leq M_2(k) \leq \omega\left(b,\frac{\eta_k}C\right).
\end{align*}
By Proposition~\ref{sne-q}, we get that
$$d(x_m,T_mx_m) \leq \frac{\eta_k}C,$$
so
$$d(x_m,T_0x_m) \leq \left( 2+\frac{\gamma_0}{\gamma_{m}} \right) d(x_m,T_mx_m) \leq \left( 2+\frac{\wt{\gamma}_0}{\gamma} \right)\cdot \frac{\eta_k}C = \eta_k,$$
and since $k \geq c = \left\lceil\frac{64b^2}{\eps^2}\right\rceil$, we have, by Lemma~\ref{l5q}, that
$$\prv{uy_k}{x_{m}y_k} \leq \frac{\eps^2}{64}.$$

We also know that $d(y_k,y_{k'}) \leq \hat{\eps} = \frac{\eps^2}{128b}$, so
\begin{align*}
\prv{uy_{k'}}{x_{m}y_{k'}} &= \prv{uy_k}{x_{m}y_k} + \prv{uy_k}{y_ky_{k'}} + \prv{y_ky_{k'}}{x_{m}y_k'}\\
&\leq \frac{\eps^2}{64} + \frac{\eps^2}{128} + \frac{\eps^2}{128} = \frac{\eps^2}{32}.
\end{align*}
So, we have shown that, for all $m \in [n,n+\hat{g}(n)+1]$, $\prv{uy_{k'}}{x_{m}y_{k'}}\leq \frac{\eps^2}{32}$.

Using Lemma~\ref{ineq} and Lemma~\ref{l2b}, we get that, for all $i \in \N$,
\begin{align*}
d^2(x_{i+1},y_{k'}) &\leq (1-\alpha_i)^2d^2(T_ix_i,y_{k'}) + 2\alpha_{i} \prv{uy_{k'}}{x_{i+1}y_{k'}}\\
&\leq (1-\alpha_{i})^2d^2(T_ix_i,T_iy_{k'}) + 2bd(y_{k'},T_iy_{k'}) + 2\alpha_{i} \prv{uy_{k'}}{x_{i+1}y_{k'}}\\
&\leq (1-\alpha_{i})^2d^2(x_i,y_{k'}) + 2bd(y_{k'},T_iy_{k'}) + 2\alpha_{i} \prv{uy_{k'}}{x_{i+1}y_{k'}}.
\end{align*}

We now seek to apply Lemma~\ref{l10} with $\eps \mapsto \frac{\eps^2}4$, $b \mapsto b^2$, $P \mapsto K(k)$ and, for all $i$, $a_i \mapsto d(x_i,y_{k'})$, $\gamma_i \mapsto 2bd(T_iy_{k'},y_{k'})$ and $\beta_i \mapsto 2\prv{uy_{k'}}{x_{i+1}y_{k'}}$.

Note that
$$n+\hat{g}(n) = n+ g^M\left(n+S\left(\frac{\eps^2}{16b^2},n\right)+1\right)+S\left(\frac{\eps^2}{16b^2},n\right) ,$$
$$\vp\left(\frac{\eps^2}4,S,K(k),b^2\right) = K(k) + S\left(\frac{\eps^2}{16b^2},K(k)\right) +1,$$
and
$$\vp\left(\frac{\eps^2}4,S,K(k),b^2\right) + g^M\left(\vp\left(\frac{\eps^2}4,S,K(k),b^2\right)\right) = \hat{K}(k),$$
so
\begin{align*}
\sum_{i=0}^{\vp\left(\frac{\eps^2}4,S,K(k),b^2\right) + g^M\left(\vp\left(\frac{\eps^2}4,S,K(k),b^2\right)\right)} 2bd(T_iy_{k'},y_{k'}) &\leq (\hat{K}(k)+1) \cdot 2b \cdot \max_{i \leq \hat{K}(k)} d(T_iy_{k'},y_{k'})\\
&\leq  (\hat{K}(k)+1) \cdot 2b \cdot M_2(k) \\
&\leq  (\hat{K}(k)+1) \cdot 2b \cdot \frac{\eps^2}{16b(\hat{K}(k)+1)} = \frac{\eps^2}8. \\
\end{align*}

Now we may apply Lemma~\ref{l10} and we get that there is an $N\leq K(k)+ S\left(\frac{\eps^2}{16b^2},K(k)\right) +1$ such that for all $i \in [N,N+g(N)]$, $d^2(x_i,y_{k'})\leq\frac{\eps^2}4$, i.e. $d(x_i,y_{k'})\leq\frac\eps2$. Now, for all $i \in [N,N+g(N)]$, since $N\leq K(k)+ S\left(\frac{\eps^2}{16b^2},K(k)\right) +1$, and so,
$$g(N) \leq g^M\left(K(k)+ S\left(\frac{\eps^2}{16b^2},K(k)\right) +1\right),$$
we have that
$$i \leq N + g(N) \leq K(k)+ S\left(\frac{\eps^2}{16b^2},K(k)\right) +1 + g^M\left(K(k)+ S\left(\frac{\eps^2}{16b^2},K(k)\right) +1\right) = \hat{K}(k),$$
so
\begin{align*}
d(y_{k'},T_iy_{k'}) &\leq \left( 2+\frac{\gamma_{i}}{\gamma_0} \right)\cdot d(y_{k'},T_0y_{k'}) \leq \left( 2+\frac{\wt{\gamma}^M_{ \hat{K}(k)}}{\gamma} \right)\cdot \frac{b}{k'} \\
&\leq \left( 2+\frac{\wt{\gamma}^M_{ \hat{K}(k)}}{\gamma} \right) \cdot b \cdot \frac1{\left\lceil \frac{\left( 2+\frac{\wt{\gamma}^M_{ \hat{K}(k)}}{\gamma} \right) \cdot b}{M_2(k)}\right\rceil}  \leq M_2(k) \leq \frac\eps2.
\end{align*}

It remains to be shown that $N \leq \Phi$. Since we have shown before that $K(k) \leq K^*$ we have that, since $S$ is nondecreasing in the second argument,
$$S\left(\frac{\eps^2}{16b^2},K(k)\right) \leq S\left(\frac{\eps^2}{16b^2},K^*\right),$$
so
$$N \leq K(k)+ S\left(\frac{\eps^2}{16b^2},K(k)\right) +1 \leq K^*+ S\left(\frac{\eps^2}{16b^2},K^*\right) +1 = \Phi.$$
The proof is now finished.
\end{proof}

As remarked before, results -- including quantitative ones -- concerning Tikhonov-regularized algorithms may be obtained from the corresponding Halpern ones, as per \cite{LeuNic17}; see also \cite[Section 3.3]{DinPinXX} for examples which specifically concern metastability. (The study of Tikhonov-regularized algorithms was also studied from the viewpoint of proof mining in \cite{DinPin21,CheLeuXX,CheKohLeuXX}.) We may now, thus, state the corresponding quantitative version of Corollary~\ref{main-t}.

\begin{corollary}\label{main-q-t}
Define, for any $g:\N\to\N$, the function $h_g:\N\to\N$, for any $n$, by $h_g(n):=g(n+1)$, and for any $R:(0,\infty)\times \N \to \N$, the function $S_R:(0,\infty)\times \N \to \N$, for any $\eps>0$ and $m \in \N$, by $S_R(\eps,m):=R(\eps,m+1)$. Also put, for any $b$, $\gamma>0$, $(\wt{\gamma}_n) \subseteq (0,\infty)$, $(\wt{\beta}_n) \subseteq (0,1]$, $\zeta:(0,\infty) \to \N$, $R:(0,\infty)\times \N \to \N$, $\eps>0$ and $g:\N\to\N$,
$$\Theta_{b,\gamma,(\wt{\gamma}_n),(\wt{\beta}_n),\zeta,R}(\eps,g) := \Phi_{b,\gamma,(\wt{\gamma}_n),(\wt{\beta}_{n+1}),\zeta,S_R}\left(\frac\eps2, h_g\right) + 1.$$

Let $(T_n)$ be a family of self-mappings of $X$, $(\gamma_n) \subseteq (0,\infty)$ and $\gamma>0$ be such that $(T_n)$ is jointly $(P_2)$ with respect to $(\gamma_n)$ and for all $n$, $\gamma_n \geq \gamma$. Let $(\wt{\gamma}_n) \subseteq (0,\infty)$ be such that for all $n$, $\wt{\gamma}_n \geq \gamma_n$. We denote by $F$ the common fixed point set of the family $(T_n)$. Let $(\beta_n) \subseteq (0,1]$, $u\in X$ and $(y_n)\subseteq X$ be such that for all $n$,
$$y_{n+1}=T_n(\beta_n u + (1-\beta_n)y_n).$$
Let $(\wt{\beta}_n) \subseteq (0,1]$ be such that for all $n$, $\wt{\beta}_n \leq \beta_n$. Let $b \in \N^*$ and $p \in F$ be such that $2d(y_0,p) \leq b$ and $2d(u,p)\leq b$. Let $\zeta:(0,\infty) \to \N$ be such that for all $\beta>0$ and all $m \geq \zeta(\beta)$, $\beta_m\leq\beta$. Let $R:(0,\infty)\times \N \to \N$ be nondecreasing in the second argument such that for all $\eps>0$ and $m \in \N$,
$$\prod_{k=m}^{R(\eps,m)} (1-\beta_k) \leq\eps.$$
Then $\Theta_{b,\gamma,(\wt{\gamma}_n),(\wt{\beta}_n),\zeta,R}$ is a rate of metastability for $(y_n)$, i.e. for all $\eps>0$ and $g: \N \to \N$, there is an $N \leq \Theta_{b,\gamma,(\wt{\gamma}_n),(\wt{\beta}_n),\zeta,R}(\eps,g)$ such that for all $i$, $j \in [N,N+g(N)]$, $d(y_i,y_j)\leq\eps$.
\end{corollary}

\begin{proof}
For all $n$, put $x_n:=\beta_n u + (1-\beta_n) y_n$, $\alpha_n:=\beta_{n+1}$ and $\wt{\alpha}_n:=\wt{\beta}_{n+1}$. We see that, for all $n$, $\wt{\alpha}_n \leq \alpha_n$, $y_{n+1}=T_nx_n$ and
$$x_{n+1}=\beta_{n+1}u + (1-\beta_{n+1})y_{n+1} = \alpha_n u + (1-\alpha_n)T_nx_n.$$
We remark that for all $\beta>0$ and all $m \geq \zeta(\beta)$, $m+1\geq \zeta(\beta)$ and so $\alpha_m = \beta_{m+1} \leq \beta$. We also remark that $S_R$ is also nondecreasing in the second argument and that for all $\eps>0$ and $m \in \N$,
$$\prod_{k=m}^{S_R(\eps,m)} (1-\alpha_k) = \prod_{k=m}^{R(\eps,m+1)} (1-\beta_{k+1}) =  \prod_{k=m+1}^{R(\eps,m+1)+1} (1-\beta_k) \leq  \prod_{k=m+1}^{R(\eps,m+1)} (1-\beta_k)\leq\eps.$$
We see that, by Busemann convexity,
$$2d(x_0,p) \leq 2(\beta_n d(u,p) + (1-\beta_n) d(y_0,p)) \leq 2\left(\beta_n \cdot \frac b 2 + (1-\beta_n) \frac b 2 \right) = b.$$

Let $\eps>0$ and $g:\N\to\N$. We may now apply Theorem~\ref{main-q} to get that there is a $w \in X$ and an $M\in\N$ with $M+1 \leq \Theta_{b,\gamma,(\wt{\gamma}_n),(\wt{\beta}_n),\zeta,R}(\eps,g)$ such that for all $q \in [M,M+h_g(M)]$, $d(w,T_qw) \leq \eps/4$ and $d(x_q,w)\leq\eps/4$.

Take $N:=M+1 \leq \Theta_{b,\gamma,(\wt{\gamma}_n),(\wt{\beta}_n),\zeta,R}(\eps,g)$. For any $i \in [N,N+g(N)]$, we have that $i-1 \in [N-1,N-1+g(N)] = [M,M+h_g(M)]$, so
\begin{align*}
d(y_i,w) &= d(T_{i-1}x_{i-1},w)\\
&\leq d(T_{i-1}x_{i-1},T_{i-1}w) + d(T_{i-1}w,w) \\
&\leq d(x_{i-1},w) + d(T_{i-1}w,w) \leq \frac\eps4+\frac\eps4 = \frac\eps2.
\end{align*}
Thus, for any $i$, $j \in [N,N+g(N)]$,
$$d(y_i,y_j) \leq d(y_i,w) + d(y_j,w) \leq \frac\eps2+\frac\eps2 = \eps.$$
The proof is now finished.
\end{proof}

Finally, Suzuki has shown in \cite{Suz07} that this sort of convergence theorems for Halpern iterations -- even for families of mappings like in our case -- directly yield convergence theorems for the corresponding viscosity iterations; this has been recently analyzed quantitatively by Kohlenbach and Pinto \cite{KohPinXX}, and the results of that paper -- specifically Lemma 3.4, Remark 3.5 and Theorem 3.11 -- may be used to immediately derive from our results rates of metastability for the viscosity proximal point algorithm, thus further illustrating the modularity of proof mining approaches.

\section{Acknowledgements}

I would like to thank Ulrich Kohlenbach and Lauren\c tiu Leu\c stean for their suggestions.

This work has been supported by a grant of the Romanian Ministry of Research, Innovation and Digitization, CNCS/CCCDI -- UEFISCDI, project number PN-III-P1-1.1-PD-2019-0396, within PNCDI III.

\end{document}